\documentclass[reprint,showpacs,amsmath,amssymb,aps,pre,superscriptaddress,floatfix,longbibliography]{revtex4-2}
\usepackage{amsthm}
\usepackage{amsmath, mathtools}
\usepackage[utf8]{inputenc}	
\usepackage[english]{babel}
\usepackage{amsmath}
\usepackage{graphicx}		
\usepackage{natbib}
\usepackage{textcomp}
\usepackage{gensymb}
\usepackage[usenames,dvipsnames,svgnames,table]{xcolor}
\usepackage[hidelinks,colorlinks=false,urlcolor=Cerulean,citecolor=black]{hyperref}
\usepackage{siunitx}
\usepackage{mathrsfs}
\usepackage{multirow}
\usepackage{bm}
\usepackage{xfrac}
\usepackage{comment}
\usepackage[normalem]{ulem}

\renewcommand{\i}{\mathrm{i}}

\DeclareMathOperator{\R}{\mathbb{R}}
\DeclareMathOperator{\C}{\mathbb{C}}

\begin{document}

\title{Geometric perspective of linear stability in finite networks of nonlinear oscillators}

\author{Yashee Sinha}
\affiliation{Department of Mathematics, Western University, London ON, Canada}
\affiliation{Birla Institute of Technology and Science, Pilani, Goa, India}

\author{Priya B. Jain}
\affiliation{Department of Mathematics, Western University, London ON, Canada}
\affiliation{Western Institute for Neuroscience, Western University, London ON, Canada}
\affiliation{Fields Lab for Network Science, Fields Institute, Toronto ON, Canada}

\author{Antonio Mihara}
\affiliation{Department of Physics, Universidade Federal de São Paulo, São Paulo, Brazil}

\author{Rene O. Medrano-T}
\affiliation{Department of Physics, Universidade Federal de São Paulo, São Paulo, Brazil}
\affiliation{Department of Physics, Universidade Estadual Paulista, São Paulo, Brazil}

\author{Ján Mináč}
\affiliation{Department of Mathematics, Western University, London ON, Canada}
\affiliation{Fields Lab for Network Science, Fields Institute, Toronto ON, Canada}

\author{Lyle E. Muller}
\thanks{These authors jointly supervised this work (co-senior authors)}
\affiliation{Department of Mathematics, Western University, London ON, Canada}
\affiliation{Western Institute for Neuroscience, Western University, London ON, Canada}
\affiliation{Fields Lab for Network Science, Fields Institute, Toronto ON, Canada}

\author{Roberto C. Budzinski}
\thanks{These authors jointly supervised this work (co-senior authors)}
\affiliation{Department of Mathematics, Western University, London ON, Canada}
\affiliation{Western Institute for Neuroscience, Western University, London ON, Canada}
\affiliation{Fields Lab for Network Science, Fields Institute, Toronto ON, Canada}

\begin{abstract}
We use a complex-valued transformation of the Kuramoto model to develop an operator-description of the linear stability in finite networks of nonlinear oscillators. This mathematical approach offers analytical predictions for the linear stability of $q$-states, which include phase synchronization ($q = 0$) and waves with different spatial frequencies ($|q| > 0$). This approach seamlessly incorporates the presence of time delays (represented by phase-lags in the coupling). With this, we are able to analytically determine the specific combination of connectivity and time delays (phase-lags) that leads to any given $q$-state to be linearly stable. This approach offers a geometric perspective of linear stability in finite networks in terms of the connectivity and delays (phase-lag), and it opens a path to designing and controlling the spatiotemporal dynamics of individual oscillator networks.
\end{abstract} 

\maketitle

The emergence of organized and sophisticated spatiotemporal dynamics is observed in several systems in nature, spanning from population of fireflies \cite{buck1976synchronous} and neural systems \cite{fell2011role,muller2018cortical} to transmission networks \cite{motter2013spontaneous}, coupled lasers \cite{nixon2011synchronized}, and nanoelectromechanical oscillators \cite{matheny2019exotic}. In this context, early works have established mathematical models to explore this behavior, where the Kuramoto model plays a central role in theoretical and numerical studies of synchronization in networks of coupled oscillators \cite{acebron2005kuramoto,rodrigues2016kuramoto,strogatz2000kuramoto}. 

In a given network of Kuramoto oscillators many different solutions can coexist, which is known as \textit{multistability} \cite{delabays2017size,wiley2006size}. In these networks, not only phase synchronization, but several phase-locking states \cite{townsend2020dense,zhang2021basins} and even exotic and more sophisticated states are possible solutions \cite{canale2015exotic}. Not all these solutions are stable, however, and studying their linear stability is an important problem that has received great attention in the past years \cite{dorfler2014synchronization,delabays2017size,wiley2006size,townsend2020dense,mehta2015algebraic,zhang2021basins,taylor2012there,jadbabaie2004stability,sokolov2019sync,mihara2022sparsity}. Despite significant progress in this issue, important questions still remain open. For instance, in a finite network of Kuramoto oscillators, what is the minimum number of connections to ensure that phase synchronization is the only stable solution \cite{townsend2020dense,kassabov2021sufficiently,lu2020synchronization}? Further, in many natural systems, the presence of delayed interactions is inherent, which is known to contribute to the emergence of spatiotemporal dynamics \cite{jeong2002time,ko2004wave,budzinski2023analytical,petkoski2022normalizing}.

Important theoretical insight on the linear stability in Kuramoto systems has been obtained using the continuum limit assumption \cite{wiley2006size,delabays2017size}. In this case, instead of a discrete set of nodes in the network, one can consider a distribution, which allows for analytical description the dynamics. This development, however, is not guaranteed to hold for finite systems. Importantly, new attempts have recently opened the possibility for an analytical framework to study the linear stability of finite Kuramoto networks without the assumption of the continuum limit \cite{mihara2019stability}. With this idea, for instance, it is possible to relate the local properties of the eigenvalues of the Jacobian matrix near equilibria with global properties like basin size \cite{mihara2022basin}. Moreover, phase-lag and time delays in the interactions between oscillators have been recently shown to change the linear stability of equilibria in Kuramoto networks \cite{mihara2019stability, an2024stability, lee2024stability}.

Here, we introduce a new approach for determining the linear stability of phase-synchronized states and wave solutions in finite networks of Kuramoto oscillators. This approach reveals how the combination of network connectivity and time delays (or phase-lag) affect the stability of these states in finite networks. Specifically, our approach is based in a complex-valued transformation to the Kuramoto model that has been recently shown to offer an operator-description of the spatiotemporal dynamics of Kuramoto networks \cite{budzinski2022geometry,budzinski2023analytical}. Using this framework, we express the linear stability of all $q$-states in terms of the eigenvalues of a composite matrix that combines the connectivity (adjacency matrix) and delays. By explicitly considering the discrete nature of finite networks, this approach offers new insights into the combinations of connectivity and delays that lead to the stability of different $q$-states.

\subsection*{Linear stability of $q$-states in finite networks}

Here, we study finite Kuramoto networks, which are described by
\begin{equation}
    \dot{\theta}_{j}(t) = \epsilon \sum_{k=0}^{N-1} A_{jk} \sin{\big(\theta_{k}(t) - \theta_{j}(t)  \big)},
    \label{eq:main_kuramoto}
\end{equation}
where $\theta_{j}(t)$ is the phase of the $j^{\mathrm{th}}$ oscillator at time $t$, $N$ is the number of nodes, $\epsilon$ is the coupling strength, and $A_{jk}$ represents the elements of the adjacency matrix. In the case where the connectivity rule is the same across nodes (circulant adjacency matrix), the solutions $\bm{\theta}^{(q)} = \big(0, \frac{2\pi q}{N}, \cdots, \frac{2\pi q(N-1)}{N} \big)$ or $q$-states are equilibria for the Kuramoto system described in Eq.~(\ref{eq:main_kuramoto}) \cite{taylor2012there,townsend2020dense,wiley2006size,nguyen2023equilibria}. Specifically, $q = 0$ represents the phase-synchronized state, where all oscillators have the same phase and frequency, and the other solutions represent traveling waves or $q$-twisted states, where the oscillators have a constant phase offset across nodes. The solutions for a given $q^{\ast}$ state appear in pairs with same spatial frequency but opposite direction of phase offset ($q = q^{\ast}$ and $q = -q^{\ast}$).

Although all of these solutions are equilibria, they can be linearly stable or unstable. The standard way to study the stability of these states is to consider a small perturbation $\eta$ at a given equilibrium point and then study how this perturbation evolves over time. The perturbed solution is given by: $\theta_{j}^{(q)} = \frac{2\pi q}{N}j + \eta_{j}$, which we then use in Eq.~(\ref{eq:main_kuramoto}) leading to $\dot{\eta}_{j}(t) = \epsilon \sum_{k = 0}^{N-1} A_{jk} \sin{\big(2\pi q N^{-1} (k-j) + \eta_{k} - \eta_{j} \big)}$. We can re-write this equation using specific trigonometric identities and the fact that $\eta \ll 1$ to obtain a linearized version given by $\dot{\eta}_j = \epsilon\sum_{k=0}^{N-1} A_{jk} \cos{(2\pi qN^{-1}(k-j))}[\eta_k - \eta_j]$ (see Sec.~\ref{sec:linear_stability_appendix}). We can now write a solution for $\eta$ in the form $\eta_{j}(t) = e^{\lambda t}e^{\frac{2\pi \i m}{N} j}$, where $\i = \sqrt{-1}$, and $m \in [0, N-1]$. With this, if $\lambda > 0$, $\eta_{j}(t)$ increases over time, and thus the solution $\bm{\theta}^{(q)}$ is unstable.

If the adjacency matrix is given by a circulant matrix, we can define $\lambda$ in terms of the generating vector of the adjacency matrix $\bm{h}$ (see  Sec.~\ref{sec:linear_stability_appendix}):
\begin{equation}
    \lambda_{m,q} = \dfrac{\hat{H}(q+m) + \hat{H}(q-m)}2 - \hat{H}(q),
    \label{eq:eigenvalues_stability}
\end{equation}
where $m$ represents the mode number of the perturbation, $q$ is the integer associated with the solution $\bm{\theta}^{(q)}$ whose stability we are investigating, and $\hat{H}(q) = \sum_{j =0}^{N-1} h_{j} e^{\frac{2\pi \i q}{N} j}$ is the discrete Fourier transform. With this, $\lambda_{m,q}$ determines the linear stability of the solution $\theta^{(q)}$ in the direction $m$. This equation has been explored in the \textit{continuum} limit, where $N \rightarrow \infty$, and $\hat{H}$ is given by the Fourier transform \cite{wiley2006size}, it is also consistent with the previous results for discrete networks \cite{mihara2019stability}. If $\lambda_{m,q} < 0$ for all $m \in [1,N-1]$, the solution $\theta^{(q)}$ is linearly stable; if $\lambda_{m,q} > 0$ in one or more directions $m$, the solution $\theta^{(q)}$ is unstable. We note that $\lambda_{0,q} = 0$ and the linear stability of the $q$-state can be determined by looking at the directions $m \in [1,N-1]$ \cite{mihara2019stability}.

\begin{figure}[tbh]
    \centering
    \includegraphics[width=0.95\columnwidth]{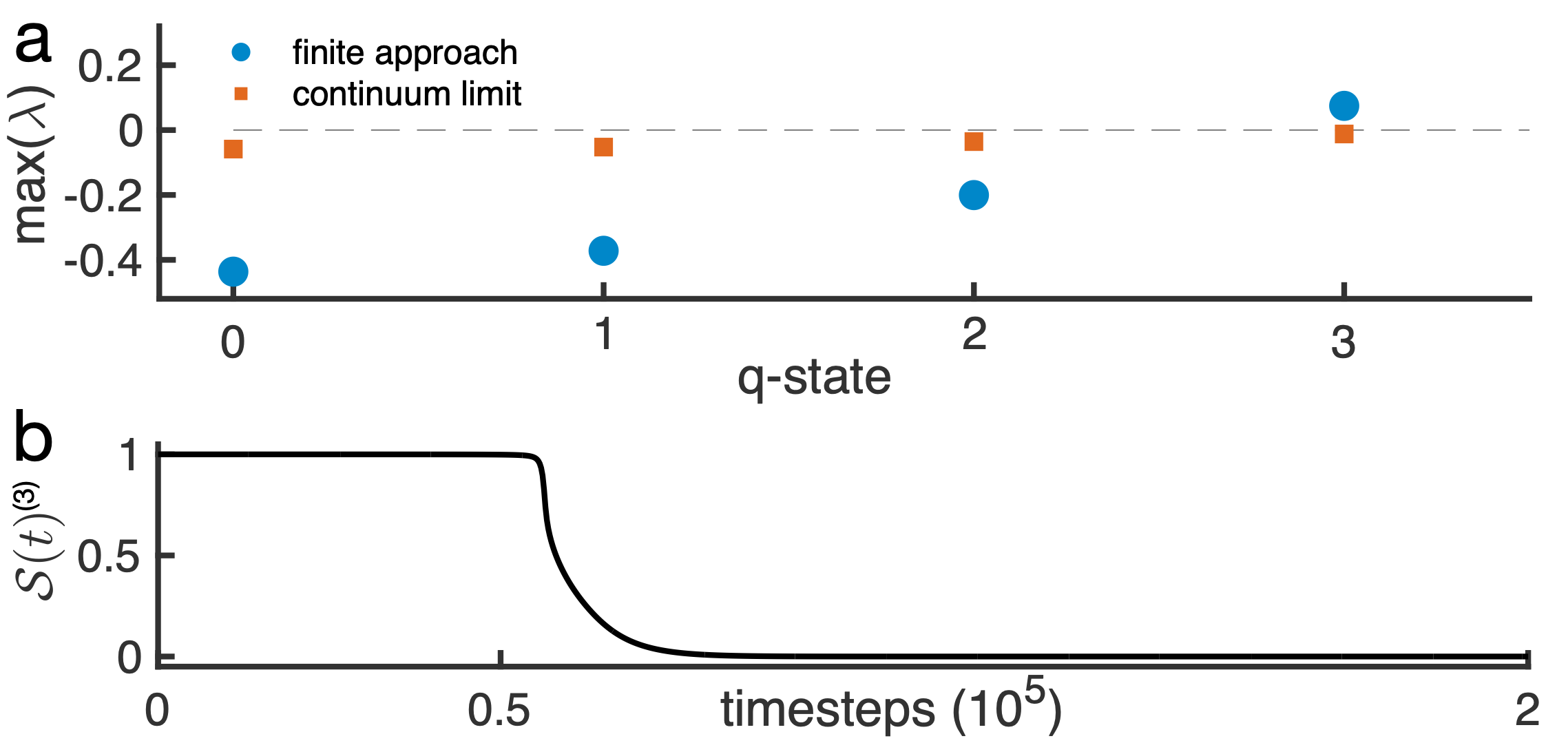}
    \caption{\textbf{Linear stability of $q$-states in finite networks.} \textbf{(a)} We consider a network with $N = 21$ oscillators on a k-ring graph with k = 2. We first use our finite approach given by Eq.~(\ref{eq:eigenvalues_stability}) and plot the maximum eigenvalue ($\mathrm{max}(\lambda)$) for each $q$-state. We observe that $q=3$-state is the first state to be unstable (blue circles). This is in contrast with the continuum approach, which predicts $q=3$-state to be linearly stable (orange square). \textbf{(b)} We then numerically integrate Eq.~(\ref{eq:main_kuramoto}) and use, as initial condition, $\bm{\theta}^{(3)}$ in addition to a small perturbation. We observe that, after $50,000$ timesteps, the network transitions to a different state, as the similarity measurement $\mathcal{S}(t)^{(3)}$ that varies from $1$ to $0$.}
    \label{fig:stability_finite_continuum}
\end{figure}

We then consider a simple example of a network defined on a k-ring graph, where each node is connected to k neighbors on each direction of the ring where the network is defined. We use our finite approach given by Eq.~(\ref{eq:eigenvalues_stability}) and plot the maximum eigenvalue $\mathrm{max}(\lambda)$ over the $N-1$ directions $m \in [1, N-1]$ for each $q$-state (Fig.~\ref{fig:stability_finite_continuum}a, blue circles). Our results predict that, for this specific network, the $q \in [0,2]$ states are linearly stable and that $q = 3$ becomes unstable. When we use the continuum limit approach as defined in \cite{wiley2006size}, where Eq.~(\ref{eq:eigenvalues_stability}) can be used with $\hat{H}(q) = \frac{\sin{(2\pi q \mathrm{k} N^{-1})}}{2\pi q \mathrm{k} N^{-1}}$, we observe that the prediction indicates the $q = 3$ state to be linearly stable (Fig.~\ref{fig:stability_finite_continuum}a, orange squares). 

To numerically test the linear stability of this equilibrium point, we integrate Eq.~(\ref{eq:main_kuramoto}) using, as initial condition, $\bm{\theta}^{(3)}$ in addition to a small perturbation ($\eta \sim 10^{-3}$). To quantify if the network remains at the initial state or transitions to a different one, we use a similarity measurement 
\begin{equation}
    \mathcal{S}(t)^{(q)} = \frac{1}{N}\left|\sum_{j=0}^{N-1} \exp{\Big(\i \theta_{j}(t)- \i\theta_{j}^{(q)}\Big)}\right|
    \label{eq:similarity}
\end{equation}
between the network state $\bm{\theta}(t)$ at time $t$ and the solution $\bm{\theta}^{(q)}$, where $\mathcal{S}^{(q)} = 1$ indicates the network is at the state $\bm{\theta}^{(q)}$ at time $t$. In the case studied in Fig.~\ref{fig:stability_finite_continuum}, $\mathcal{S}(t)^{(3)}$ transitions from $1$ to $0$ (Fig.~\ref{fig:stability_finite_continuum}b). These numerical results confirm that the $q=3$ state is, in fact, unstable, which agrees with our prediction given by the finite approach.

\subsection*{Operator description of linear stability in finite networks}

In this work, we use an operator-based description of Kuramoto networks to show that the linear stability of any $q$-state in any circulant network can be directly determined by the eigenvalues of an aggregated matrix that combines connectivity and phase-lag (or time delays). We first introduce the operator-based description of Kuramoto networks in the case without phase-lag or delays, and then extend to these cases in later sections of the paper.

We recently introduced a complex-valued transformation to the Kuramoto model, which is able to reproduce the trajectories of finite Kuramoto networks under a variety of dynamical behaviors, including phase synchronization, traveling waves (or phase-locking states), and even chimera states \cite{budzinski2022geometry,budzinski2023analytical}. This approach can be described as an iterative process given by the combination of two operators:
\begin{equation}
    \bm{x}(t+\varsigma) = \Lambda \big[e^{\varsigma \bm{K}} \bm{x}(t)\big],
    \label{eq:cv_approach}
\end{equation}
where $\bm{x} \in \C^{N}$, $\varsigma$ is small but finite, fixed value, $\Lambda: x \rightarrow \sfrac{x}{|x|}$ represents an elementwise operator mapping the modulus of each state vector element $x_j$ to unity, and $\bm{K} = \epsilon \bm{A}$, with $N$, $\epsilon$, and $\bm{A}$ being the same as described in the Kuramoto model. Importantly, the trajectories of the Kuramoto network given by $\theta_{j}(t)$ can be compared to $\mathrm{Arg}[x_{j}(t)]$, and these two systems are precisely equivalent when $\sfrac{x_{k}}{x_{j}} = 1$ \cite{budzinski2022geometry}. Equation (\ref{eq:cv_approach}) reproduces the spatiotemporal dynamics of Kuramoto networks using a combination of a linear operator, the matrix exponential, and a nonlinear operator, the operator $\Lambda$. Importantly, the change in the phase dynamics of $\bm{x}(t)$ is given by the matrix exponential, so that we can analyze the nonlinear dynamics of finite Kuramoto networks directly in terms of eigenspectrum of $\bm{K}$ \cite{budzinski2023analytical, budzinski2022geometry}. In the case the adjacency matrix $\bm{A}$ is circulant (and so is $\bm{K}$), the eigenvalues and eigenvectors are determined analytically by the circulant diagonalization theorem \cite{davis1979}. The $j^{\mathrm{th}}$ eigenvalue is given by $\varphi_{j} = \sum_{k=0}^{N-1} h_{k} \exp{\big(\frac{-2\pi \i}{N} jk \big)}$ and the $s$ entry of the $j^{\mathrm{th}}$ eigenvector is given by $(\bm{v}_{j})_s = (\sfrac{1}{\sqrt{N}}) \exp{\big( \frac{-2\pi \i}{N} j s \big)}$. In particular, the argument of these eigenvectors precisely represents the solutions $\bm{\theta}^{(q)}$ \cite{nguyen2023equilibria}. Moreover, the eigenvalues, $\varphi_j$, are equal to $\hat{H}(q)$ defined in Eq.~(\ref{eq:eigenvalues_stability}). Here, we show that linear stability of $q$-states in circulant networks of Kuramoto oscillators can be directly interpreted and determined by analyzing the eigenvalues of $\bm{K}$
\begin{equation}
   \lambda_{m,q} = \frac{\gamma_{q+m} + \gamma_{q-m}}{2} - \gamma_{q},
    \label{eq:eigenvalues_cdt_stability}
\end{equation}
where $\gamma$ are the real part of the eigenvalues of $\bm{K}$ ($\gamma_{j} = \mathrm{Real}[\varphi_{j}]$) and we consider the ordering to be $N$-periodic, i.e. $q\pm m \equiv (q\pm m)~\mathrm{mod}~N$. Importantly, Eq.~(\ref{eq:eigenvalues_cdt_stability}) is equivalent to Eq.~(\ref{eq:eigenvalues_stability}) and so is its interpretation of the linear stability of a given equilibrium point given by $\mathrm{Arg}[\bm{v}_{q}] = \bm{\theta}^{(q)}$: if $\lambda_{m,q} \leq 0, \, \forall \, m \in [1,N-1]$, then $\bm{\theta}^{(q)}$ is linearly stable, and $\bm{\theta}^{(q)}$ is unstable otherwise. Moreover, here we show that the largest eigenvalue of $\bm{K}$ in real part determines the $q$-state with largest basin of attraction. Importantly, with this new formulation, we are able to seamlessly consider the case of heterogeneous phase-lag or even distance-dependent delays (as shown in Secs.~\ref{sec:phase_lag} and \ref{sec:delays}). With this, the matrix $\bm{K}$ can aggregate connectivity and phase-lag or delays, and by using Eq.~(\ref{eq:eigenvalues_stability}) we are able to determine the linear stability of all $q$-state in the presence of phase-lags or time delays and the solution with largest basin of attraction. Thus, this offers a new, fully analytical way to describe the linear stability for the equilibria in any finite Kuramoto model on any circulant network just by looking at the connectivity pattern and phase-lag or time delays in the interactions in the network.

\subsection*{Linear stability in k-ring graphs}

\begin{figure}[t!]
    \centering
    \includegraphics[width=\columnwidth]{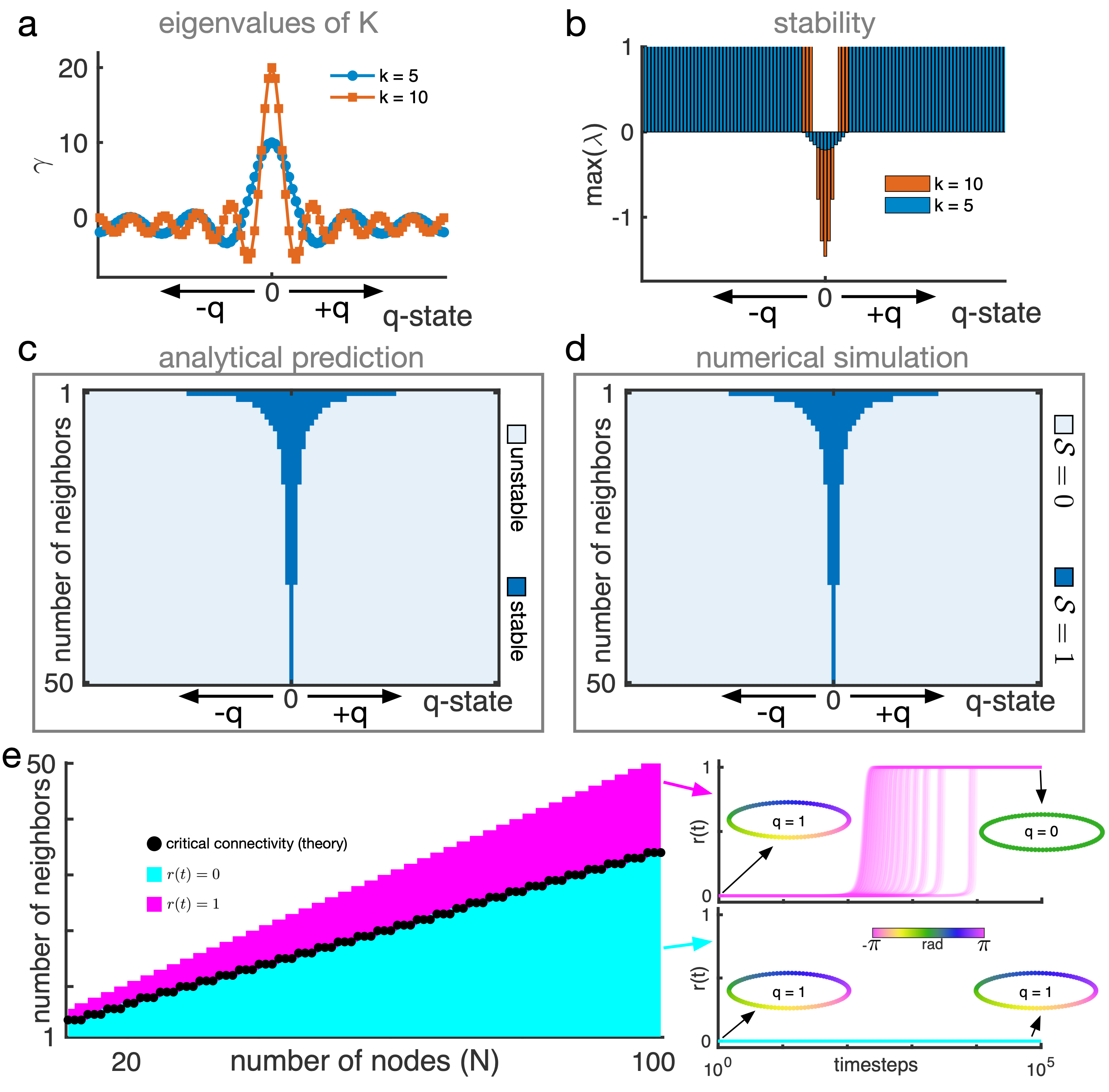}
    \caption{\textbf{Network connectivity and stability.} \textbf{(a)} We consider k-ring networks and calculate the eigenvalues of the matrix $\bm{K}$. \textbf{(b)} We then use Eq.~(\ref{eq:eigenvalues_cdt_stability}) to determine the linear stability of the $q$-states based on the eigenvalues $\gamma$. We can observe how the changes in $\gamma$ due to the connectivity change the stability properties. \textbf{(c)} We consider the full range of $k$ from the first neighbors k = 1 to the complete graph (all-to-all) k = 50, where dark blue indicates the linearly stable states and light blues indicate unstable states. \textbf{(d)} We numerically integrate Eq.~(\ref{eq:main_kuramoto}) using the $q$-state $\bm{\theta}^{(q)}$ in addition to a small perturbation as initial state and then measure the similarity measurement $\mathcal{S}^{(q)}$ at the end of the simulation at $t = 2 \times 10^{5}$. \textbf{(e)} We can analytically determine the  connectivity (represented by $k$) to ensure that phase synchronization is the only linearly stable solution among all the $q$-states (black circles, left). To test our predictions, we integrate Eq.~(\ref{eq:main_kuramoto}) using $\bm{\theta}^{(1)}$ in addition to a small perturbation as the initial state. We then plot the order parameter $r(t)$ in the end of the simulation (color-code, left), which shows that our predictions are correct. We also show examples of the time evolution of the system $r(t)$ for different conditions (right).}
    \label{fig:stability_eigenvalues_cdt}
\end{figure}

As an example, we consider k-ring networks, where k is the number of neighbors that a given oscillator is connected to in each direction on the ring. Importantly, by changing k we change the number of connections in the network and also the eigenvalues of $\bm{K}$ (Fig.~\ref{fig:stability_eigenvalues_cdt}a), but the eigenvectors remain the same, where the first eigenvector represents phase synchronization ($q = 0$-state) and the others represent different $q$-states or twisted states. We then use Eq.~(\ref{eq:eigenvalues_cdt_stability}) to analytically predict the linear stability of these equilibrium points. We observe that the changes in the connectivity (different k values) lead to changes in the eigenvalues $\gamma$, which then change the stability of the $q$-states (Fig.~\ref{fig:stability_eigenvalues_cdt}b). We then generate analytical predictions for the stability of all $q$-states for networks going from first neighbors (k = 1) to global networks, where all oscillators are connected (k = 50). In Fig.~\ref{fig:stability_eigenvalues_cdt}c, dark blue represents linearly stable states and light blue represents unstable ones. We observe that, as k increases, the $q$-states lose stability and, after a critical $k$ value, only the phase synchronized state ($q = 0$) remains linearly stable. 

In order to test our analytical predictions, we numerically integrate Eq.~(\ref{eq:main_kuramoto}) using the $q$ solutions $\bm{\theta}^{(q)}$ in addition to a small perturbation as the initial condition for the network on the k-ring graphs with k varying from k = 1 to k = 50. We integrate the system for $t = 2 \times 10^{5}$ timesteps and evaluate the similarity measurement $\mathcal{S}^{(q)}$ at the end of the simulation. In Fig.~\ref{fig:stability_eigenvalues_cdt}d, dark blue represents $\mathcal{S}^{(q)} = 1$, which indicates the $q$-state is linearly stable, since the network returns to the $q$-state after the perturbation; and light blue represents $\mathcal{S}^{(q)} = 0$, which indicates the $q$-state is unstable, since the network transitions to a different state after the perturbation. These results show that our analytical predictions match the numerical simulations.

Figures \ref{fig:stability_eigenvalues_cdt}c and \ref{fig:stability_eigenvalues_cdt}d show that, for low values of $k$, where the network has few connections, many $q$-states are linearly stable; but as the number of connections increases these $q$-states lose stability and after a specific value of $k$, only the phase synchronized $q = 0$ state is linearly stable. This result leads to an important question: what are the necessary conditions to ensure that phase synchronization will emerge in a given network? Our results suggest that, for finite networks, the critical connectivity where only the $q = 0$ state is stable is discrete and varies discretely with the size of the network $N$ for k-ring graphs. 

As an example, we consider the k-ring networks with different number of nodes $N \in [10, 100]$ and vary the degree of the network k (note that k $ \in [1, \mathrm{floor}(\sfrac{N}{2})]$). Figure \ref{fig:stability_eigenvalues_cdt}e shows our analytical predictions (black circles) for the critical k when the $q = 0$ state is the only stable state among the $q$-states, which shows discrete values that resemble a staircase. We then numerically tested our predictions by integrating Eq.~(\ref{eq:main_kuramoto}) under the same conditions and considering the initial state given by $\bm{\theta}^{(1)}$ in addition to a small perturbation. We then measure the Kuramoto order parameter $r(t) = \frac{1}{N}\Big|\sum_{j=0}^{N-1} e^{\i \theta_{j}(t)}\Big|$ and consider the final value of $r(t)$ after $t = 1 \times 10^{5}$ timesteps. Figure \ref{fig:stability_eigenvalues_cdt}e shows, in color-code, the order parameter for the final state in the simulation, where $r = 0$ indicates the network remains at $\bm{\theta}^{(1)}$ and this state is still linearly stable, and $r = 1$ indicates the network transitions to phase synchronization. These results show that, the connectivity to ensure that phase synchronization is the only linearly stable solution among all the $q$-states varies discretely with network size in finite k-ring graphs.

\subsection*{Linear stability in networks with distance-dependent coupling}

\begin{figure}[tbh]
    \centering
    \includegraphics[width=\columnwidth]{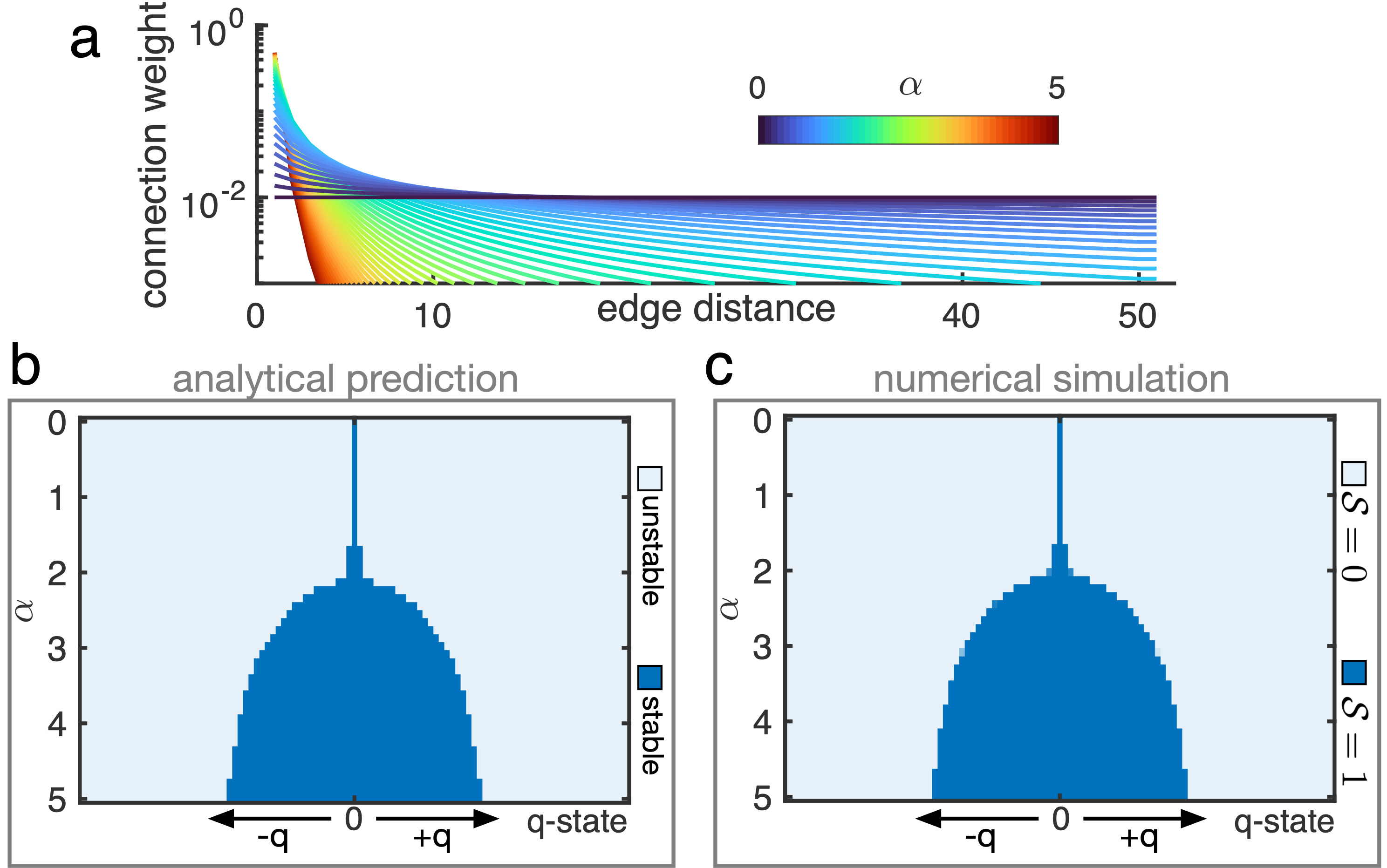}
    \caption{\textbf{Linear stability of finite networks on weighted graphs.} \textbf{(a)} We consider distance-dependent networks, where the connection weight between two nodes decay with their edge distance. The parameter $\alpha$ controls this decay: for $\alpha = 0$, all nodes are connected with the same weight; as $\alpha$ increases, near nodes are connected with stronger weights. \textbf{(b)} We use Eq.~(\ref{eq:eigenvalues_cdt_stability}) to analytically predict the linear stability of $q$-states for distance-dependent networks with different $\alpha$ values. \textbf{(c)} We observe that our prediction are correct by using numerical simulations for these networks and evaluating the similarity measurements $\mathcal{S}^{(q)}$ at the end of the simulation ($t = 2.5 \times 10^5$).}
    \label{fig:distance_dependent}
\end{figure}

We also apply our approach to study weighted adjacency matrices. Here, we consider a distance-dependent connectivity, where the nodes are arranged on an one dimensional ring with periodic boundary conditions, and the weight of the connection between two nodes $j$ and $k$ decay with the edge distance $d_{jk} = \mathrm{min}(|j - k|, N - |j - k|)$. The adjacency matrix is given by $A_{jk} =\frac{1}{\rho(\alpha) (d_{jk})^\alpha}$, where $\alpha$ is the exponent that controls the decay of the connection weights and $\rho$ is a normalization factor such that the row sum of $\bm{A}$ is always equal to one. Figure \ref{fig:distance_dependent}a shows the connection weights as a function of the edge distance $d$ for different $\alpha$ values (in color-code). We observe that, for $\alpha = 0$, the all nodes are connected with the same weight (global network), but as $\alpha$ increases, the nodes display stronger connections weight near neighbors. 

Because this network is circulant, we can directly apply our approach -- given by Eq.~(\ref{eq:eigenvalues_cdt_stability}) -- to analytically obtain the linear stability for the $q$-states in these networks (Fig.~\ref{fig:distance_dependent}b). We observe that, for low values of $\alpha$, only the phase synchronized state is linearly stable, but as $\alpha$ increases, and the network becomes more locally connected, different $q$-states become linearly stable. We then perform a numerical simulation for these network and verify the stability of these states, where we find a great match with out analytical predictions (Fig.~\ref{fig:distance_dependent}c). These results we obtain here are similar to the results found in \cite{lee2024stability} where the linear stability is calculated through the eigenvalues of the Jacobian of the system.

\subsection*{Linear stability in networks with heterogeneous phase-lag}\label{sec:phase_lag}

One of the main advantages of the operator-based description of the Kuramoto model is the fact that the complex-valued transformation easily incorporates the presence of phase-lag (homogeneous or heterogeneous) in the interaction term \cite{budzinski2022geometry,budzinski2023analytical}. In this case, the Kuramoto model, or Sakaguchi-Kuramoto model, can be written as:
\begin{equation}
    \dot{\theta}_{j}(t) = \epsilon \sum_{k=0}^{N-1} A_{jk} \sin{\big(\theta_{k}(t) - \theta_{j}(t) -\phi_{jk} \big)},
    \label{eq:kuramoto_phase_lag}
\end{equation}
where $\phi_{jk}$ is the phase-lag in the interaction between oscillators $j$ and $k$.

The complex-valued transformation is given by Eq.~(\ref{eq:cv_approach}), where the matrix $\bm{K}$ incorporates the phase-lag effect on the systems. Specifically, we have that $K_{jk} = \epsilon e^{-\i\phi_{jk}} A_{jk}$ \cite{budzinski2023analytical}. Importantly, if the matrix representing the phase-lag $\phi$ is circulant (if the rule is the same across node -- homogeneous, distance dependent, and others) and the adjacency matrix $\bm{A}$ is circulant, the matrix $\bm{K}$ is also circulant. This means that under the presence of phase-lag the matrices $\bm{A}$ and $\bm{K}$ share the same eigenvectors, and thus the $q$-states remain as solutions for Eq.~(\ref{eq:kuramoto_phase_lag}). The presence of phase-lag in the coupling changes, however, the eigenvalues of $\bm{K}$, which then changes the linear stability of the $q$-states under the presence of phase-lag. Our approach given by Eq.~(\ref{eq:eigenvalues_cdt_stability}) is able to determine the linear stability of these states under the presence or absence of phase-lag. Specifically, the phase-lag $\phi$ will rotate the eigenvalues of $\bm{K}$ in the complex plane, and our approach is able to analytically and geometrically show how exactly phase-lag is changing the linear stability of $q$-states in Sakaguchi-Kuramoto networks.

\begin{figure}[t]
    \centering
    \includegraphics[width=\columnwidth]{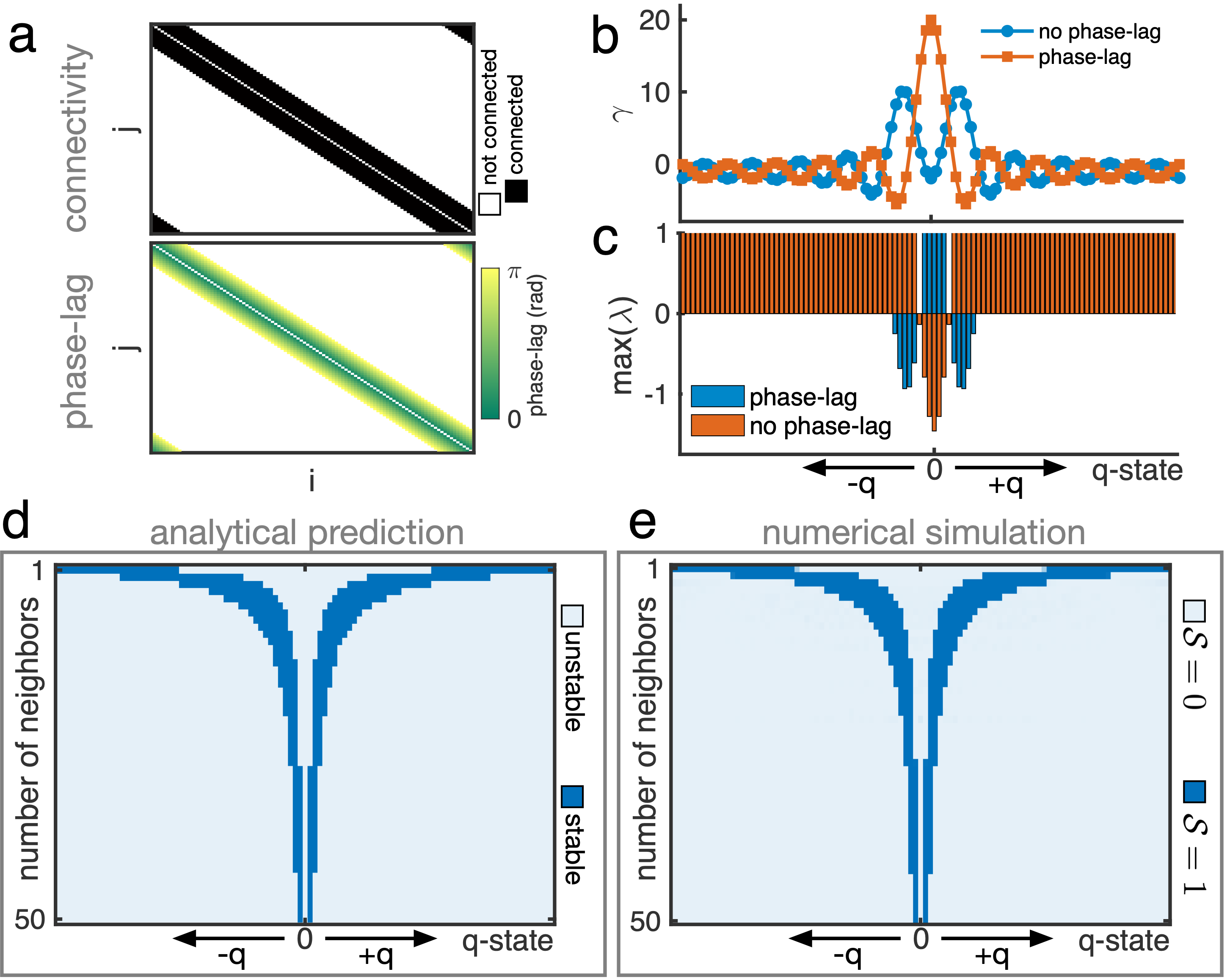}
    \caption{\textbf{Analytical predictions for the linear stability of $q$-states in networks with phase-lags.} \textbf{(a)} Here, we consider k-ring networks, where the pattern of connections is binary (top) and phase-lag is given by $\phi_{jk} = \frac{\pi d_{jk}}{\mathrm{k}}$ (bottom), which varies in $[0,\,\pi]$ (color-code) and no phase-lag is considered between two nodes that are not connected (represented in white). \textbf{(b)} As a specific example, we consider the case with $N = 101$ and k = 10, with and without phase-lag in the coupling, and calculate the eigenvalues $\gamma$ of the matrix $\bm{K}$. \textbf{(c)} We then use Eq.~(\ref{eq:eigenvalues_cdt_stability}) to determine the linear stability of the $q$-states in the case with phase-lags (blue bars) and without phase-lags (orange bars), where we observe different stable $q$-states due to the presence of phase-lags. We then consider the full range of k $\in [1, 50]$ for these networks with heterogeneous phase-lag, \textbf{(d)} where our analytical predictions given by Eq.~(\ref{eq:eigenvalues_cdt_stability}) match the \textbf{(e)} numerical simulations of Eq.~(\ref{eq:kuramoto_phase_lag}) for testing the stability of the $q$-states.}
    \label{fig:stability_phase_lag}
\end{figure}

We consider the same networks on k-ring graphs as studied in Fig.~\ref{fig:stability_eigenvalues_cdt} but with heterogeneous phase-lag. Here, the adjacency matrix is binary (Fig.~\ref{fig:stability_phase_lag}a, top), and the phase-lag matrix is given by a distance-dependent rule such that $\phi_{jk} = \frac{\pi d_{jk}}{\mathrm{k}}$, where $d_{jk} = \mathrm{min}(|j - k|, N - |j - k|)$ (Fig.~\ref{fig:stability_phase_lag}a, bottom). The presence of heterogeneous phase-lag only affects the eigenvalues $\varphi$, by inducing rotations in the complex plane, which then lead to changes in the real part of the eigenvalues $\gamma$ (Fig.~\ref{fig:stability_phase_lag}b). This, in turn, affects the linear stability of $q$-states (Fig.~\ref{fig:stability_phase_lag}c). For this specific case ($N = 101$, k$ = 10$), without phase-lag the $q$-states with $q = 0, q = \pm 1, q = \pm 2, q = \pm 3$. For the same network, when we consider the presence of distance-dependent phase-lag, we observe that the linearly stable $q$-states become $q = \pm 4, q = \pm 5, q = \pm 6, q = \pm 7, q = \pm 8$.

We then considered different networks where we vary the number of neighbors k. By using Eq.~(\ref{eq:eigenvalues_cdt_stability}), we are able to analytically predict the linear stability of the $q$-states under the presence of heterogeneous phase-lag (Fig.~\ref{fig:stability_phase_lag}d). We observe for low k, high-spatial frequency $q$-states are linearly stable, and as k increases, they lose stability and lower $q$-states become linearly stable. However, the phase synchronized state is not linearly stable at any value of k (in contrast to the case without phase-lag depicted in Fig.~\ref{fig:stability_eigenvalues_cdt}). To test these predictions, we numerically integrate Eq.~(\ref{eq:kuramoto_phase_lag}) using as initial conditions the $\bm{\theta}^{(q)}$ states in addition to a small perturbation. We then integrate the system for $t = 2 \times 10^{5}$ and evaluate the similarity measurement $\mathcal{S}^{(q)}$ at the end of the simulation (Fig.~\ref{fig:stability_phase_lag}e), which shows the agreement between the analytical predictions and the numerical results. This also shows the specific conditions with respect for k and $\bm{\phi}$ in which any given $q$-state becomes the only linearly stable state among all the $q$-states.

\subsection*{Linear stability in networks with heterogeneous delay}\label{sec:delays}

As a direct application of our analytical approach, we can now predict the spatial frequency of traveling waves that emerge in finite networks with multiple time delays, which is represented by
\begin{equation}
    \theta_{j}(t) = \omega + \sum_{k=0}^{N-1} A_{jk}\sin{\big(\theta_{k}(t - \tau_{jk}) - \theta_{j}(t) \big)}.
    \label{eq:kuramoto_delays}
\end{equation}
Time delays can be approximated by phase-lags, where $\theta_{k}(t - \tau_{jk}) \approx \theta_k(t) - \omega\tau_{jk}$, with $\omega$ being the natural frequency of oscillation. In this case, the matrix $\bm{K}$ can be written as $K_{jk} = \epsilon\exp{(-\i\omega\tau_{jk})}A_{jk}$, where we consider the time delays to linearly increase with the distance between nodes such that $\tau_{jk} = \sfrac{d_{jk}}{\nu}$, with $d_{jk} = \mathrm{min}(|j - k|, N - |j - k|)$ and $\nu$ representing the conduction speed of the signal. Because the matrix representing the time delays is also circulant in this case, we can write the real part of the eigenvalues of $\bm{K}$ explicitly as (see Sec.~\ref{sec:stability_delays_appendix} for details): 
\begin{equation}
   \gamma_j = \sum_{k=0}^{N} \epsilon A_{0k} \cos{\bigg(\frac{-2\pi}{N} jk - \omega\tau_{0k} \bigg)}. 
\end{equation}
This can be directly applied to Eq.~(\ref{eq:eigenvalues_cdt_stability}) to analytical determine the linear stability of the $q$-state solutions explicitly as a function of connectivity ($\bm{A}$) and time delays ($\bm{\tau}$). 
\begin{figure}[thb]
    \centering
    \includegraphics[width=\columnwidth]{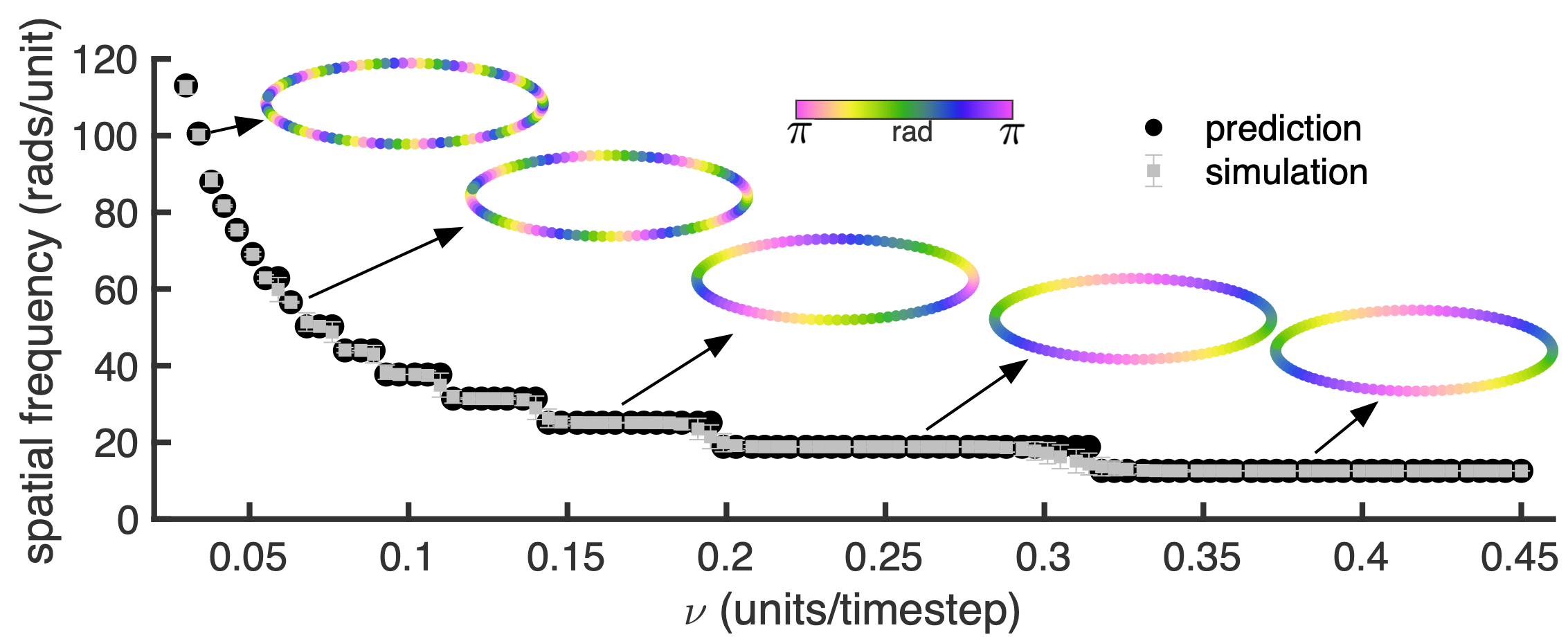}
    \caption{\textbf{Analytical predictions for the emergence of waves in delayed networks.} We consider a k-ring graph with $N = 101$ and k = 50 (global network), $\epsilon = 1$, and different values of conduction speed $\nu$. We use Eq.~(\ref{eq:eigenvalues_cdt_stability}) to analytically predict the stability of the $q$-states in the delayed system as the conduction speed is varied. The largest eigenvalue $\gamma$ allows us to estimate the stable $q$-state with largest basin of attraction, which gives us a way to predict the spatial frequency of the wave that emerges for a given $\nu$ (black dots). We then numerically integrate the delayed Kuramoto model Eq.~(\ref{eq:kuramoto_delays}) with random initial conditions and measure the spatial frequency of the dynamics after $t = 1 \times 10^{5}$ timesteps (gray squares). We repeat this procedure for $1,000$ different initial states, where the gray bars represent the standard deviation.}
    \label{fig:waves_speed}
\end{figure}

We then use this to study the emergence of waves in delayed system as the conduction speed $\nu$ is varied, which affects the times delays in the network. Based on Eq.~(\ref{eq:eigenvalues_cdt_stability}) we obtain the linear stability for the possible solutions for the system. Using the ideas introduced in \cite{mihara2022basin}, we can directly use $\gamma$ to estimate which stable $q$-state has the largest basin size. With this, we can now predict the spatial frequency of the waves that emerge in these delayed systems for different conduction speeds by simply looking at the largest eigenvalue $\gamma$ (Fig.~\ref{fig:waves_speed}, black dots). In this example, we observe that, the spatial frequency of the wave decreases as the conduction speed increases where the discreteness of the system plays an important role (staircase-like behavior for higher conduction speeds). To test our predictions, we numerically integrate the Kuramoto network with time delays -- Eq.~(\ref{eq:kuramoto_delays}) -- starting with a random initial state and then analyze the spatiotemporal pattern after $t = 1 \times 10^{5}$ timesteps (gray squares, Fig.~\ref{fig:waves_speed}), which show great agreement with our predictions.

\section*{Discussion}

In this work, we have introduced a new geometric and analytical view on the linear stability of $q$-states in finite networks of Kuramoto oscillators on circulant adjacency matrices. This approach offers a clear and simple way to analytically predict the linear stability of any $q$-state in finite networks just by looking at the eigenvalues of the aggregate matrix $\bm{K}$, which combines the connectivity of the network and the presence of phase-lags or time delays in the interactions. This approach is possible due to our operator-based description of Kuramoto networks \cite{budzinski2022geometry,budzinski2023analytical}, which links directly the emergent nonlinear dynamics of nonlinear oscillator networks to the eigenspectrum of $\bm{K}$. Importantly, this matrix $\bm{K}$ is composed of the adjacency matrix $\bm{A}$, the coupling strength $\epsilon$ and it seamless incorporates the presence of (heterogeneous) phase-lag $\bm{\phi}$, which is a good approximation for networks with delayed interactions. Further, our approach also naturally generalizes to networks on weighted graphs. With this, we are able to look at the combination of connectivity and phase-lags and analytically predict the linear stability of equilibrium points in finite networks of Kuramoto oscillators with delayed interactions. There has been recent interest on this topic with important results regarding the linear stability \cite{lee2024stability,an2024stability,mihara2019stability,mihara2022basin} and the critical connectivity to ensure phase synchronization \cite{townsend2020dense,yoneda2021lower,kassabov2021sufficiently}, our framework offers a new perspective for this problem. 

Particularly, this framework offers a path to analytically determine the combination of connectivity and delays (phase-lags) so that phase synchronization or a wave solution ($q$-states) is the only linearly stable solution among all the possible $q$-states for any finite Kuramoto network on any circulant matrix. This results may offer a different view on the critical connectivity problem, where we embrace the discrete nature of finite networks and obtain a discrete description of the role of single connections being added or having their weight changed in a given network system. The results in this work shed light on an important problem in network theory and nonlinear dynamics: how does a pattern of connections in combination with nonlinear coupling lead to the emergence of specific spatiotemporal dynamics. This question is of great importance in neuroscience \cite{curto2019relating,carandini2012circuits} with potential applications in computations and machine learning \cite{masoliver2022embedded,budzinski2024exact,bottcher2022ai,menara2022functional}. Because our approach considers finite networks and seamlessly incorporates time delays (through the phase-lag approximation), it helps to create a new perspective on this problem where standard techniques from linear algebra can lead to important results for this field.

\section{Appendix}

\subsection{Linear stability analysis for finite networks}\label{sec:linear_stability_appendix}

Here, we study the linear stability of the $q$-states, given by $\bm{\theta}^{(q)} = \big(0, \frac{2\pi q}{N}, \cdots, \frac{2\pi q(N-1)}{N} \big)$. We consider a small perturbation $\bm{\eta} \in \R^N$ with $|\bm{\eta}|\ll 1$ around the twisted states such that $\bm{\theta}^{'(q)} = \bm{\theta}^{(q)} + \bm{\eta}$. If we consider this perturbed solution in Eq.~(\ref{eq:main_kuramoto}), we obtain:
\begin{equation}
    \dot{\eta}_j = \epsilon \sum_{k=0}^{N-1} A_{jk} \sin{ \bigg(\dfrac{2\pi q}N(k-j) + \eta_k - \eta_j \bigg)},
\end{equation}
\begin{equation}
    \begin{split}
        \dot{\eta}_j = \epsilon \sum_{k=0}^{N-1} A_{jk} \Bigg[\sin\bigg(\dfrac{2\pi q}N(k-j)\bigg)\cos(\eta_k - \eta_j ) \\ + \cos\bigg(\dfrac{2\pi q}N(k-j)\bigg)\sin(\eta_k - \eta_j)\Bigg].
    \end{split}
\label{eq:perturbation}
\end{equation}
Because $\eta_j$ are very small, we can consider $\sin(\eta_j) \approx \eta_j$ and $\cos(\eta_j) \approx 1$. Furthermore, because $\bm{\theta}^{(q)}$ is an equilibrium point,  $\sum_{k=0}^{N-1} A_{jk} \sin{\big(2\pi qN^{-1}(k-j)\big)} = 0$. With this, we obtain
\begin{equation}
    \dot{\eta}_j = \epsilon\sum_{k=0}^{N-1} A_{jk} \cos{\bigg(\dfrac{2\pi q}N(k-j)\bigg)} \big[\eta_k - \eta_j \big].
    \label{eq:linear_eta}
\end{equation}
The analytical approach we develop is valid for any circulant network, which is the case where the connectivity rule is the same across nodes \cite{davis1979}, where the adjacency matrix is circulant. In this case, the adjacency matrix $\bm{A}$ is composed by cyclic permutation of the generating vector $\bm{h} = (h_0, h_1, h_2, \ldots, h_{N-1})$, such that $\bm{A} = \mathrm{circ}(\bm{h})$.

We first define $\bm{M} \in \R^{N\times N}$ as $M_{kj} = A_{jk}\cos\big(2\pi qN^{-1} (k-j)\big)$ and $\bm{1} = [1,1,\dots, 1] \in \R^N$. Using the $\text{diag}(\cdot)$ operator, which produces a diagonal matrix with the elements of a given vector, we reformulate the Eq.~(\ref{eq:linear_eta}) in its matrix form:
\begin{align}
    \dot{\bm{\eta}} &= \epsilon \big(\bm{M}\bm{\eta} - \text{diag}(\bm{\eta})\bm{M}\bm{1}\big)\label{eq:matrix}
\end{align}
Since $\bm{M}$ is also circulant, all rows of the matrix are made of the same elements and therefore the sum of the elements in every row is identical. Therefore, $\bm{M}\bm{1} = \sum_{j=0}^{N-1}M_{0j} \bm{1}$. Putting this in Eq.~(\ref{eq:matrix}) we get,
\begin{align}
    \dot{\bm{\eta}} &= \epsilon\left(\bm{M} - \sum_{j=0}^{N-1}M_{0j}\right)\bm{\eta}.\label{eq:preeigen}
\end{align}

To find the eigenvectors of the RHS, we utilize the circulant diagonalization theorem (CDT) \cite{davis1979}. According to CDT, the eigenvalues of a circulant matrix, $\bm{A} = \text{circ}(\bm{h})$, are given by the discrete Fourier transform of the generating vector, $\bm{h}$:
\begin{equation}
    \hat{H}(q)=  \sum\limits_{k=0}^{N-1} h_k \exp\left[ - \frac{2\pi \i}{N} jk \right],
\label{eq:ktheigval}
\end{equation}
and the eigenvectors are given by:
\begin{equation}
    [\bm{v_m}]_j = \exp\left(\dfrac{-2\pi\i}N mj\right).
\end{equation}
Therefore, the eigenvalues of the matrix $\bm{M}$ are given by: 
\begin{equation}
    \sum_{j=0}^{N-1} M_{0j} e^{\frac{\i2\pi m}N j} = \sum_{j=0}^{N-1} h_{j} \cos\left(\frac{2\pi q}N j\right)\exp\left(\frac{\i2\pi m}N j\right)
\end{equation}
or
\begin{equation}
    \sum_{j=0}^{N-1} M_{0j} e^{\frac{\i2\pi m}N j} = \frac{\hat{H}(q+m) + \hat{H}(-q + m)}2.
\end{equation}
If we assume $\bm{\eta} = f(t)\bm{v}_m$ for some $m\in \R$, and $f: \R \rightarrow\R$, Eq.~(\ref{eq:preeigen}) becomes
\begin{equation}
    \dot{\bm{\eta}} = \epsilon\left(\frac{\hat{H}(q+m) + \hat{H}(-q + m)}2 -  \sum_{j=0}^{N-1} M_{0j}\right)\bm{\eta} .
\end{equation}
The solutions for this set of equation is then given by:
\begin{equation}
    \bm{\eta} = \exp(\lambda t) \bm{v}_m,
\end{equation}
where 
\begin{equation}
    \lambda = \frac{\hat{H}(q+m) + \hat{H}(-q + m)}2 -  \sum_{j=0}^{N-1} M_{0j}.
\end{equation}

The form of $\lambda$ can be further simplified by noting that since $\bm{h}\in \R^N$ and is also symmetric around $\lfloor N/2 \rfloor$ ($h_i = h_{N-i}$), the DFT $\hat{H}(m)$ is an even function around zero, i.e. $\hat{H}(m) = \hat{H}(-m)$ \cite{Wongbook}. We also note that 
\begin{equation}
    \sum_{j=0}^{N-1}M_{0j} = \sum_{j=0}^{N-1}h_j\cos(2\pi q N^{-1} j),
\end{equation}
and
\begin{equation}
   \sum_{j=0}^{N-1}h_j\cos(2\pi q N^{-1} j) = \frac{\hat{H}(q)+\hat{H}(-q)}2 = \hat{H}(q).
\end{equation}

Therefore, the eigenvalues can be written as:
\begin{align}
    \lambda_{m,q} = \frac{\hat{H}(q+m) + \hat{H}(q - m)}2 -  \hat{H}(q).
\end{align}
For any given $q$-state, its linear stability depends on whether the eigenvalues, $\lambda_{m,q}$, are negative for all $m \in \{1,2, \dots N\}$. 

\subsection{Analytical description of the eigenspectrum in the presence of multiple time delays}\label{sec:stability_delays_appendix}

To determine the linear stability of $q$-states in finite networks of Kuramoto oscillators with distance-dependent delays, we consider the complex-valued transformation described in Eq.~(\ref{eq:cv_approach}), where the matrix $\bm{K}$ is given by $\bm{K} = \epsilon e^{{(-\i\omega\tau_{jk})}} A_{jk}$, where $\bm{\tau}$ represents the delay matrix and $\bm{A}$ the adjacency matrix. If the matrices $\bm{A}$ and $\bm{\tau}$ are circulants, the matrix $\bm{K}$ is also circulant. In this case, the generating vector is represented by $h_{j} = e^{(-\i\omega\tau_{0j})}A_{0j}$. With this, we can write the eigenvalues of $\bm{K}$ as
\begin{equation}
    \varphi_{j} =  \sum\limits_{k=0}^{N-1} e^{(-\i\omega\tau_{0k})}A_{0k} e^{ - \frac{2\pi \i}{N} jk }.
\end{equation}
From this, we can combine the exponential terms, which leads to
\begin{equation}
    \varphi_{j} =  \sum\limits_{k=0}^{N-1} A_{0k} \exp\left[ -\i \Bigg( \frac{2\pi}{N} jk + \omega\tau_{0k} \Bigg) \right].
\end{equation}
The eigenvalue with largest real part allows us to predict the spatiotemporal pattern that emerges in the network, which is given by the argument of the eigenvector associated to this eigenvalue \cite{budzinski2023analytical}. With this in mind, we can derive an expression for the real part of the eigenvalues of $\bm{K}$ as a function of the time delays:
\begin{equation}
    \gamma_j = \sum\limits_{k=0}^{N-1} A_{0k} \cos{\Bigg(\frac{-2\pi}{N} jk - \omega\tau_{0k} \Bigg).}
\end{equation}
This equation allows us to analytically express the largest eigenvalue (in real part) for a given set of heterogeneous time delays, thus allowing us to predict the specific wave number, or spatial frequency of the wave that will emerge in the dynamics. Since, $\mathrm{Real}[\lambda_{m,q}]$ describes the stability of the $q$-states and the phase of the eigenvector $\bm{v}_q$ represents the $q$-state,
\begin{equation}
    \lambda_{m,q} = \frac{\gamma_{(q+m)} + \gamma_{(q-m)}}{2} - \gamma_q
\end{equation}    
describes the stability of the $q$-state for the Kuramoto networks with distance-dependent delays.

\vspace{0.75cm}

\textit{Acknowledgments:} This work was supported by BrainsCAN at Western University through the Canada First Research Excellence Fund (CFREF), the NSF through a NeuroNex award (\#2015276), the Natural Sciences and Engineering Research Council of Canada (NSERC) grant R0370A01, Compute Ontario (computeontario.ca), and Digital Research Alliance of Canada (alliancecan.ca). A.M. acknowledges the financial support of the São Paulo Research Foundation (FAPESP), grant number \#2023/08144-3. R.O.M acknowledges the support of National Council for Scientific and Technological Development – CNPq, project number 408522/2023-2.

\end{document}